\documentclass[12pt]{article}

\usepackage{amssymb,amsmath}
\usepackage{url}

\usepackage{color}

\DeclareMathAlphabet{\eufrak}{U}{}{}{} 
\SetMathAlphabet\eufrak{normal}{U}{euf}{m}{n}
\SetMathAlphabet\eufrak{bold}{U}{euf}{b}{n}

\numberwithin{equation}{section}

\newenvironment{Proof}{\removelastskip\par\medskip
\noindent{\em Proof.} \rm}{\penalty-20\null\hfill$\square$\par\medbreak}

 \allowdisplaybreaks

 \def\real{{\mathord{\mathbb R}}}
 
 \def\inte{{\mathord{\mathbb N}}}
 
 \def\qu{{\mathord{\mathbb Z}}}
  
 \def\Cov{{\mathrm{{\rm Cov \kern1pt}}}}

 \def\real{{\mathord{{\rm I\kern-3pt R}}}}        
 
 \def\inte{{\mathord{{\rm I\kern-3pt N}}}}
 \def\sZZ{{\rm Z\kern-.45em{}Z}}

 \def\sQQ{{\kern 0.27em \vrule height1.45ex width0.03em depth0em
           \kern-0.30em \rm Q}}
 \def\qu{{\mathchoice
         {\sQQ}
         {\sQQ}
   {\kern 0.225em \vrule height1.05ex width0.025em depth0em \kern-0.25em \rm Q}
   {\kern 0.180em \vrule height0.78ex width0.020em depth0em \kern-0.20em \rm Q}
         }}
 \def\sGG{{\kern 0.27em \vrule height1.45ex width0.03em depth0em
           \kern-0.30em \rm G}}
 \def\gg{{\mathchoice
         {\sGG}
         {\sGG}
   {\kern 0.225em \vrule height1.05ex width0.025em depth0em \kern-0.25em \rm G}
   {\kern 0.180em \vrule height0.78ex width0.020em depth0em \kern-0.20em \rm G}
         }}

 \newtheorem{prop}{Proposition}[section]

 \newtheorem{theorem}[prop]{Theorem}

\def\E{\mathop{\hbox{\rm I\kern-0.20em E}}\nolimits}

 \newcounter{hyp}
\title{\huge 
 Mixing of Poisson random measures under interacting transformations 
} 
 
\date{\normalsize \today 
}

\author{
Nicolas Privault\thanks{nprivault@ntu.edu.sg}
\\ 
\normalsize 
Division of Mathematical Sciences 
\\ 
\normalsize 
School of Physical and Mathematical Sciences 
\\ 
\normalsize 
Nanyang Technological University 
\\ 
\normalsize 
SPMS-MAS, 21 Nanyang Link 
\\ 
\normalsize 
Singapore 637371
}

\begin{document}

\maketitle 

\vspace{-0.8cm} 
 
\begin{abstract}
{ 
We derive sufficient conditions for the mixing of all orders of interacting transformations of a spatial Poisson point process, under a zero-type condition in probability and a generalized adaptedness condition. This extends a classical result in the case of deterministic transformations of Poisson measures. The approach relies on moment and covariance identities for Poisson stochastic integrals with random integrands. 
} 
\end{abstract}
\noindent {\bf Key words:} Poisson random measures; interacting transformations; mixing; ergodicity. 
\\ 
{\em Mathematics Subject Classification:} 37A25, 60G57, 37A05, 60H07. 

\baselineskip0.7cm
 
\hyphenation{func-tio-nals} 
\hyphenation{Privault} 

\section{Introduction} 
\label{s123} 
\noindent 
 The ergodicity and mixing properties of Poisson random 
 measures under deterministic transformations 
 have been considered by several authors, 
 cf. e.g. \cite{marchat}, \cite{grabinsky}, \cite{roy}. 
 This paper investigates mixing beyond the deterministic case 
 by considering interacting, i.e. configuration dependent, 
 transformations of Poisson samples. 
\noindent 
\\ 
 
 Consider a $\sigma$-compact metric space 
 $X$ with Borel $\sigma$-algebra ${\cal B}(X)$ 
 and let $\Omega$ denote the configuration space on 
 $(X,{\cal B} (X))$, i.e. 
$$ 
 \Omega = \left\{ 
 \omega = 
 ( x_i )_{i=1}^N \subset X, \ x_i\not= x_j \ 
 \forall i\not= j, \ N \in \inte \cup \{ \infty \} \right\} 
, 
$$ 
 is the space of at most countable subsets 
 of $X$, whose elements $\omega \in \Omega$ are identified to the 
 point measures 
\begin{equation} 
\label{dsp} 
 \displaystyle  
 \omega (dy) = \sum_{x \in \omega} \delta_x (dy) 
, 
\end{equation} 
 where $\delta_x$ denotes the Dirac measure at $x\in X$. 
 The space $\Omega$ is endowed with the Poisson probability 
 measure $\pi_\sigma$ with $\sigma$-finite diffuse intensity 
 $\sigma (dx)$ on $X$ 
 and its associated $\sigma$-algebra ${\cal F}$ 
 generated by $\omega \mapsto \omega (A)$ for $A\in 
 {\cal B}(X)$ such that $\sigma (A) < \infty$. 
 In particular, $\pi_\sigma (d\omega )$-almost surely, 
 $\omega \in \Omega$ is 
 locally finite on compact sets and \eqref{dsp} is a Radon measure. 
\\ 
 
 Given a measurable random transformation 
$$ 
 \tau : X \times \Omega \longrightarrow X
,
$$ 
 of $X$ and an 
 element $\omega $ of $\Omega$ of the form \eqref{dsp}, 
 let $\tau_* (\omega )$ denote the transformation of 
 $\omega \in \Omega$ by $\tau ( \cdot , \omega ) : X \longrightarrow X$, i.e. 
\begin{equation} 
\label{llk2} 
 \tau_* (\omega) 
 : = \sum_{x \in \omega } 
 \delta_{\tau ( x, \omega ) }, 
 \qquad 
 \omega \in \Omega,  
\end{equation} 
 is the image measure of $\omega (dy)$ by 
 $\tau ( \cdot , \omega ) : X \longrightarrow X$. 
 In other words, the transformation 
\begin{equation} 
\label{llk} 
 \tau_* : \Omega \longrightarrow \Omega 
\end{equation} 
 shifts every configuration point $x \in \omega$ according to 
 $x \longmapsto \tau ( x, \omega )$, and in the deterministic case 
 $\tau_*$ is also called the 
 Poisson suspension over $\tau : X \longrightarrow X$, cf. 
 \S~9.1 of \cite{cornfeld}. 
\\ 
 
 In Theorem~4.8 of \cite{roy} it is shown, 
 using the moment generating function of Poisson 
 random measures, that 
 a conservative deterministic dynamical system 
 $(\Omega, \pi_\sigma, \sigma , \tau )$ 
 where $\tau : X \longrightarrow X$ leaves $\sigma$ 
 invariant 
 is mixing of all orders if and only if 
 $\tau : X \longrightarrow X$ is of zero type, i.e. 
$$ 
 \lim_{n \to \infty} \langle h , h \circ \tau^n \rangle_{L^2_\sigma (X)} 
 = 
 0 
, 
$$ 
 for all $h\in L^2_\sigma (X)$, 
 cf. also \cite{cornfeld} for the Gaussian case. 
\\ 
 
 In Theorem~\ref{prnt} below we show that an interacting transformation 
 $\tau ( \cdot , \omega ) : X \longrightarrow X$ 
 leaving $\sigma$ invariant $\pi_\sigma (d\omega )$-a.s. 
 is mixing of all orders provided the family of transformations 
 $\tau^{(n)} : X \times \Omega \longrightarrow \ \Omega$, 
 $n \in \inte$, inductively defined by 
 $\tau^{(0)} ( x, \omega ) : =  x$ and 
\begin{equation} 
\label{rr*} 
 \tau^{(n)} ( x , \omega ) 
 : = 
 \tau^{(n-1)} ( \tau ( x , \omega ) , \tau_* \omega ) 
, 
 \qquad 
 n \geq 1, 
\end{equation} 
 $\omega \in \Omega$, $x \in X$, 
 satisfies the zero-type condition 
$$ 
 \lim_{n \to \infty} \langle 
 g 
 , 
 h \circ \tau^{(n)} 
 \rangle_{L^2_\sigma (X)} = 0 
$$ 
 in probability for all $g , h\in {\cal C}_c (X)$, 
 as well as the vanishing gradient condition \eqref{111.111.111} 
 below that plays the role of an adaptedness condition in the 
 absence of time ordering. 
\\ 
 
 When $\tau : X \longrightarrow X$ is deterministic, 
 Condition~\eqref{111.111.111} is always satisfied 
 and we have 
$$ 
 \tau^{(n)} ( x , \omega ) = \tau^n ( x ), 
 \qquad 
 \omega \in \Omega, \quad 
 x \in X, \quad 
 n \geq 1, 
$$ 
 hence Theorem~\ref{prnt} recovers the classical mixing 
 conditions on the Poisson space as 
 it suffices to state Condition~\eqref{condfg} for $g=h$, 
 in which case it becomes equivalent 
 to the deterministic zero-type condition 
$$ 
 \lim_{n \to \infty} \langle 
 h 
 , 
 h \circ \tau^n 
 \rangle_{L^2_\sigma (X)} = 0, 
 \qquad 
 h\in {\cal C}_c (X). 
$$ 
 Our proof uses extension to joint moments 
 of the moment identities for Poisson stochastic integrals 
 with random integrands of \cite{momentpoi}, 
 cf. \cite{flint} for an extension to point processes. 
\\ 
 
 Related arguments have been previously applied on the Wiener space 
 using the Skorohod integral, cf. \cite{prmix}, 
 \cite{uzergocrs}, \cite{uzergo}. 
\\ 
 
\indent 
 This paper is organized as follows. 
 In Section~\ref{s00} we state and 
 recall some preliminary results on 
 invariance of Poisson random measures 
 and joint moment identities for 
 Poisson stochastic integrals. 
 In Section~\ref{exs} we present and prove our main result on the 
 mixing property of interacting transformations. 
 In Section~\ref{examples} 
 we consider a family of examples based on 
 transformations conditioned by the random boundary of a convex Poisson hull. The invariance of such transformations with respect to the Poisson 
 measure is consistent with the intuitive fact that 
 the distribution of the inside points remains Poisson 
 when they are shifted within its convex hull 
 according to the data of the vertices, 
 cf. the unpublished manuscript \cite{davydov}. 
\section{Invariance and joint moment identities} 
\label{s00} 
 In this section we recall some preliminary results on 
 invariance of Poisson random measures under interacting transformations, 
 and we derive joint moment 
 identities for the 
 Poisson stochastic integral $\int_X u(x,\omega ) \omega (dx)$ 
 of a random integrand $u: X \times \Omega \longrightarrow \real$. 
\subsubsection*{Invariance of Poisson random measures} 
  Let now $D_x$, $x\in X$, 
 denote the finite difference gradient 
 defined for all $\omega \in \Omega$ and $x\in X$ as 
$$ 
 D_x F(\omega ) 
 = F(\omega \cup \{ x \} ) - F( \omega ) 
, 
$$ 
 for any random variable $F:\Omega \longrightarrow \real$, 
 cf. e.g. 
 Theorem~6.5 page 21 of \cite{yito}. 
 Given $\Theta = \{ x_{k_1},\ldots ,x_{k_l} \} \subset \{ x_1, \ldots , x_n \}$ 
 and $u: X^n \times \Omega \longrightarrow \real$, 
 we define the iterated gradient 
\begin{equation} 
\label{djkld1} 
 D_\Theta 
 u (x_1 , \ldots , x_n , \omega )  
 : = 
 D_{x_{k_1}} \cdots  D_{x_{k_l}}  u (x_1 , \ldots , x_n , \omega ), 
 \qquad 
 x_1,\ldots ,x_n \in X. 
\end{equation} 
 Recall that by Theorem~3.3 of \cite{prinv} or 
 \cite{prigirid}, or Theorem~5.2 of \cite{bretonprivaultfact}, 
 $\tau_* : \Omega \longrightarrow \Omega$ leaves $\pi_\sigma$ invariant, 
 i.e. 
 $\tau_* \pi_\sigma = \pi_\sigma$,
 provided 
 that for $\pi_\sigma$-a.s. $\omega \in \Omega$ 
 the random transformation 
 $\tau ( \cdot , \omega ) : X \longrightarrow X$ 
 leaves $\sigma (dx)$ invariant 
 and satisfies the vanishing condition 
\begin{equation} 
\label{cyclic4} 
 D_{\Theta_1} 
 \tau ( x_1 , \omega ) 
 \cdots 
 D_{\Theta_m} 
 \tau ( x_m , \omega ) 
 = 
 0, 
\end{equation} 
 for every family $\{ \Theta_1, \ldots , \Theta_m \}$ 
 of (non empty) subsets such that 
 $\Theta_1 \cup \cdots \cup \Theta_m = \{x_1,\ldots ,x_m \}$, 
 for all $x_1,\ldots ,x_m \in X$, 
 $\pi_\sigma ( d \omega )-a.s.$, 
 $m\geq 1$. 
\\ 
 
 Condition~\eqref{cyclic4} is interpreted 
 by saying that 
 for $\omega \in \Omega$ and 
 $x_1,\ldots ,x_m \in X$ 
 there exists $l \in \{1,\ldots , m \}$ such that 
\begin{equation} 
\label{cc}  
 D_{x_l} \tau ( x_{(l+1) \! \! \! \! \mod m} , \omega ) = 0, 
 \quad 
 i.e. 
 \quad 
 \tau ( x_{(l+1) \! \! \! \! \mod m} , \omega \cup \{ x_l \} ) 
 = 
 \tau ( x_{(l+1) \! \! \! \! \mod m} , \omega ) 
, 
\end{equation} 
 where 
$(l \! \mod m )= l$, $1\leq l \leq m$, 
 and 
$(m+1 \! \mod m )= 1$, 
 i.e. the $m$-tuples 
$$ 
 ( 
 \tau ( x_2 , \omega \cup \{ x_1\} ), 
 \tau ( x_3 , \omega \cup \{ x_2\} ), 
 \ldots 
 , 
 \tau ( x_m , \omega \cup \{ x_{m-1} \} ) 
 , 
 \tau ( x_1 , \omega \cup \{ x_m \} ) 
 ) 
$$ 
 and 
$ 
 ( 
 \tau ( x_2 , \omega ), 
 \tau ( x_3 , \omega ), 
 \ldots 
 , 
 \tau ( x_m , \omega ) 
 , 
 \tau ( x_1 , \omega ) 
 ) 
$ 
 coincide on at least one component in $X^m$, 
 cf. page 1074 of \cite{priinvcr}. When $m=1$, 
 Condition~\eqref{cyclic4} reads $D_x \tau ( x , \omega ) = 0$, 
 i.e. $\tau ( x , \omega \cup \{ x \} ) = \tau ( x , \omega )$, 
 $x\in X$, $\pi_\sigma (d\omega )$-a.s. 
\\ 
 
 Condition~\eqref{cyclic4} is known to hold 
 when $\tau : X \times \Omega \longrightarrow X$ is predictable 
 with respect to a total binary relation 
 $\preceq$ on $X$, which is the case in particular 
 when $X$ is of the form $X = \real_+ \times Z$ 
 and $\tau : X \times \Omega \longrightarrow X$ is predictable 
 with respect to the canonical filtration 
 $({\cal F}_t)_{t\in \real_+}$ 
 generated on $X = \real_+ \times Z$, 
 cf. Section~4 of \cite{prinv}. 
\subsubsection*{Joint moment identities} 
 For any random variable $F:\Omega \longrightarrow \real$, 
 we let $\varepsilon_{x_1,\ldots ,x_k}^+$ denote the addition 
 operator defined as 
$$ 
 \varepsilon_{x_1,\ldots ,x_k}^+ F(\omega ) : = 
 F(\omega \cup \{ x_1,\ldots , x_k \} )
, 
 \qquad 
 \omega \in \Omega, 
 \quad 
 x_1,\ldots ,x_k \in X. 
$$ 
 Next, given $u: X \times \Omega \longrightarrow X$ a measurable process, 
 we define the Poisson stochastic integral of $u$ as 
$$ 
 \int_X u(x,\omega ) \omega ( dx ) 
 = 
 \sum_{x\in \omega} u(x , \omega ), 
$$ 
 provided the sum converges absolutely, $\pi_\sigma (d\omega )$-a.s. 
 In the next proposition we extend 
 Proposition~3.1 of \cite{momentpoi} to a joint moment identity 
 using an induction argument. 
\begin{prop} 
\label{pr11} 
 Let $u: X \times \Omega \longrightarrow X$ be a measurable process 
 and $n=n_1+\cdots + n_p$, $p\geq 1$. 
 We have 
\begin{eqnarray} 
\label{eee} 
\lefteqn{ 
 E \left[ 
 \left( 
 \int_X u_1 (x,\omega ) \omega ( dx ) 
 \right)^{n_1} 
 \cdots 
 \left( 
 \int_X u_p (x,\omega ) \omega ( dx ) 
 \right)^{n_p} 
 \right] 
} 
\\ 
\nonumber 
 & = & 
 \sum_{k=1}^n 
 \sum_{P^n_1,\ldots , P^n_k } 
 E \left[ 
 \int_{X^k} 
 \varepsilon_{x_1,\ldots ,x_k}^+ 
 \left( 
 \prod_{j=1}^k 
 \prod_{i=1}^p 
 u_i^{l^n_{i,j}} 
 (x_j , \omega ) 
 \right)  
 \sigma (dx_1) \cdots  \sigma (dx_k) 
 \right] 
, 
\end{eqnarray} 
 where the sum runs over all partitions 
 $P^n_1,\ldots , P^n_k$ of $\{ 1 , \ldots , n \}$ 
 and the power $l^n_{i,j}$ is the cardinal 
$$ 
 l^n_{i,j} : = 
 | P^n_j \cap ( n_1+\cdots + n_{i-1},n_1+\cdots + n_i ]|,
 \qquad
 i=1,\ldots,k, \quad j=1,\ldots ,p, 
$$ 
 for any $n\geq 1$ such that all terms in the 
 right hand side of \eqref{eee} are integrable. 
\end{prop} 
\begin{Proof} 
 We will show the modified identity 
\begin{eqnarray} 
\label{eee.1} 
\lefteqn{ 
 E \left[ 
 F \left( 
 \int_X u_1 (x,\omega ) \omega ( dx ) 
 \right)^{n_1} 
 \cdots 
 \left( 
 \int_X u_p (x,\omega ) \omega ( dx ) 
 \right)^{n_p} 
 \right] 
} 
\\ 
\nonumber 
 & = & 
 \sum_{k=1}^n 
 \sum_{P^n_1,\ldots , P^n_k } 
 E \left[ 
 \int_{X^k} 
 \varepsilon_{x_1,\ldots ,x_k}^+ 
 \left( 
 F 
 \prod_{j=1}^k 
 \prod_{i=1}^p 
 u_i^{l^n_{i,j}} 
 (x_j , \omega ) 
 \right)  
 \sigma (dx_1) \cdots  \sigma (dx_k) 
 \right] 
, 
\end{eqnarray} 
 for $F$ a sufficiently integrable random variable,
 where $n=n_1+\cdots +n_p$. 
 For $p=1$ the identity is Proposition~3.1 of \cite{momentpoi}. 
 Next we assume that the identity holds at the rank 
 $p\geq 1$. 
 Replacing $F$ with $F\left( 
 \int_X u_{p+1} (x,\omega ) \omega ( dx ) 
 \right)^{n_{p+1}}$ in \eqref{eee.1} we get 
\begin{eqnarray*} 
\lefteqn{ 
 E \left[ 
 F 
 \left( 
 \int_X u_1 (x,\omega ) \omega ( dx ) 
 \right)^{n_1} 
 \cdots 
 \left( 
 \int_X u_{p+1} (x,\omega ) \omega ( dx ) 
 \right)^{n_{p+1}} 
 \right] 
} 
\\ 
\nonumber 
 & = & 
 \sum_{k=1}^n 
 \sum_{P^n_1,\ldots , P^n_k } 
 \int_{X^k} 
 \sigma (dx_1) \cdots  \sigma (dx_k) 
\\ 
\nonumber 
 & & 
 E \left[ 
 \varepsilon_{x_1,\ldots ,x_k}^+ 
 \left( 
 F 
 \left( 
 \int_X u_{p+1} (x,\omega ) \omega ( dx ) 
 \right)^{n_{p+1}} 
 \prod_{j=1}^k 
 \prod_{i=1}^p 
 u_i^{l^n_{i,j}} 
 (x_j , \omega ) 
 \right)  
 \right] 
\\ 
\nonumber 
 & = & 
 \sum_{k=1}^n 
 \sum_{P^n_1,\ldots , P^n_k } 
 \int_{X^k} 
 E \Biggl[ 
 \left( 
 \int_X 
 \varepsilon_{x_1,\ldots ,x_k}^+ 
 u_{p+1} (x,\omega ) \omega ( dx ) 
 + 
 \sum_{i=1}^k 
 \varepsilon_{x_1,\ldots ,x_k}^+ 
 u_{p+1} (x_i ,\omega ) 
 \right)^{n_{p+1}} 
\\ 
\nonumber 
 & & 
 \varepsilon_{x_1,\ldots ,x_k}^+ 
 \left( 
 F 
 \prod_{j=1}^k 
 \prod_{i=1}^p 
 u_i^{l^n_{i,j}} 
 (x_j , \omega ) 
 \right)  
 \Biggr] 
 \sigma (dx_1) \cdots  \sigma (dx_k) 
\\ 
\nonumber 
 & = & 
 \sum_{k=1}^n 
 \sum_{P^n_1,\ldots , P^n_k } 
 \sum_{a_0+\cdots +a_k=n_{p+1}} 
 \frac{n_{p+1}!}{a_0!\cdots a_k!} 
 \int_{X^k} 
 E \Biggl[ 
 \left( 
 \int_X 
 \varepsilon_{x_1,\ldots ,x_k}^+ 
 u_{p+1} (x,\omega ) \omega ( dx ) 
 \right)^{a_0} 
\\ 
\nonumber 
 & & 
 \varepsilon_{x_1,\ldots ,x_k}^+ 
 \left( 
 F 
 \prod_{j=1}^k 
 \left( 
 u_{p+1}^{a_j} (x_j ,\omega ) 
 \prod_{i=1}^p 
 u_i^{l^n_{i,j}} 
 (x_j , \omega ) 
 \right) 
 \right) 
 \Biggr] 
 \sigma (dx_1) \cdots  \sigma (dx_k) 
\\ 
\nonumber 
 & = & 
 \sum_{k=1}^n 
 \sum_{P^n_1,\ldots , P^n_k } 
 \sum_{a_0+\cdots +a_k=n_{p+1}} 
 \frac{n_{p+1}!}{a_0!\cdots a_k!} 
 \sum_{j=1}^{a_0} 
 \int_{X^{k + a_0}} 
 E \Biggl[ 
 \sum_{Q^{a_0}_j,\ldots , Q^{a_0}_j } 
\\ 
\nonumber 
 & & 
 \varepsilon_{x_1,\ldots ,x_{k+a_0}}^+ 
 \left( 
 F 
 \prod_{q=k+1}^{k+a_0} u_{p+1}^{|Q^{a_0}_q|} (x_q ,\omega ) 
 \prod_{j=1}^k 
 \left( 
 u_{p+1}^{a_j} (x_j ,\omega ) 
 \prod_{i=1}^p 
 u_i^{l^n_{i,j}} 
 (x_j , \omega ) 
 \right) 
 \right) 
 \Biggr] 
 \sigma (dx_1) \cdots  \sigma (dx_{k+a_0}) 
\\ 
\nonumber 
 & = & 
 \sum_{k=1}^{n+n_{p+1}} 
 \sum_{P^{n+n_{p+1}}_1,\ldots , P^{n+n_{p+1}}_k } 
 E \Biggl[ 
 \int_{X^k} 
 \varepsilon_{x_1,\ldots ,x_k}^+ 
 \left( 
 F 
 \prod_{l=1}^k 
 \prod_{i=1}^{p+1} 
 u_i^{l^{n+n_{p+1}}_{i,j}} 
 (x_l , \omega ) 
 \right)  
 \sigma (dx_1) \cdots  \sigma (dx_k) 
 \Biggr] 
, 
\end{eqnarray*} 
 where the summation over the partitions 
 $P^{n+n_{p+1}}_1,\ldots , P^{n+n_{p+1}}_k$ of 
 $\{ 1, \ldots , n+n_{p+1} \}$, 
 is obtained by combining the partitions 
 of $\{ 1, \ldots , n \}$ with 
 the partitions 
 $Q^{a_0}_j,\ldots , Q^{a_0}_j$ of 
 $\{ 1, \ldots , a_0 \}$ and 
 $a_1, \ldots , a_k$ elements 
 of $\{1,\ldots , n_{p+1}\}$ which are counted 
 according to $n_{p+1}! / ( a_0!\cdots a_k! )$, with 
$$ 
 l^{n+n_{p+1}}_{p+1,j} 
 = 
 l^n_{i,j} + a_j, \quad 1 \leq j \leq k, 
 \qquad 
 l^{n+n_{p+1}}_{p+1,j} 
 = 
 l^n_{i,j} + |Q^{a_0}_q|, \quad k+1 \leq j \leq k+a_0, 
. 
$$ 
\end{Proof} 
 Note that when $n=1$, \eqref{eee} coincides with 
 the classical Mecke \cite{jmecke} identity 
\begin{equation} 
\label{mecke1} 
 E 
 \left[ 
 \int_X u(x,\omega ) \omega ( dx ) 
 \right] 
 = 
 E \left[ 
 \int_{X} 
 \varepsilon_x^+ 
 u (x , \omega ) 
 \sigma (dx) 
 \right] 
. 
\end{equation} 
\section{Mixing of interacting transformations} 
\label{exs} 
 Theorem~\ref{prnt} is the main result of this paper. 
 The vanishing condition \eqref{111.111.111} below is stated 
 in the sense of \eqref{cc} above. 
\begin{theorem} 
\label{prnt} 
 Assume that 
 $\tau ( \cdot , \omega ) : X \longrightarrow X$ 
 leaves $\sigma (dx)$ invariant for $\pi_\sigma$-a.s. $\omega \in \Omega$, 
 and 
\begin{equation} 
\label{111.111.111} 
 D_{\Theta_1} 
 \tau^{(k_1)} ( x_1 , \omega ) 
 \cdots 
 D_{\Theta_m} 
 \tau^{(k_m)} ( x_m , \omega ) 
 = 
 0, 
\end{equation} 
 for every family $\{ \Theta_1, \ldots , \Theta_m \}$ 
 of (non empty) subsets such that 
 $\Theta_1 \cup \cdots \cup \Theta_m = \{ x_1,\ldots ,x_m \}$, 
 $x_1,\ldots ,x_m \in X$ and all 
 $\pi_\sigma ( d \omega )$-a.s., 
 $k_1,\ldots ,k_m \geq 1$, $m \geq 1$. 
 Then the measure-preserving transformation 
 $\tau_* : \Omega \longrightarrow \Omega$ 
 is mixing of all orders $m \geq 1$ 
 provided the zero-type condition 
\begin{equation} 
\label{condfg} 
 \lim_{n \to \infty} \langle g , h \circ \tau^{(n)} \rangle = 0 
\end{equation} 
 is satisfied in probability for all $g , h \in {\cal C}_c (X)$. 
\end{theorem} 
\begin{Proof} 
 Let $k_{i,n} := p_{1,n}+\cdots +p_{i,n}$, $i=1,\ldots ,m$, 
 where 
 $ 
 (p_{1,n})_{n \geq 1}, 
 \ldots 
 , 
 (p_{m,n})_{n \geq 1} 
$ 
 is a family of $m$ strictly increasing sequences of integers. 
 Consider $h_1,\ldots , h_m \in {\cal C}^+_c (X)$ 
 nonnegative continuous functions bounded by $1$ with compact support, 
 and let $l_1,\ldots ,l_m \geq 1$. 
 In order to prove mixing of order $m$ 
 we need to show that the joint moments 
\begin{eqnarray} 
\label{djda} 
\lefteqn{ 
 E \left[ 
 \left( 
 \int_X h_1 ( x ) \omega ( dx ) 
 \right)^{l_1} \circ \tau_*^{k_{1,n}} 
 \cdots 
 \left( 
 \int_X h_m ( x ) \omega ( dx )  
 \right)^{l_m} \circ \tau_*^{k_{m,n}} 
 \right] 
} 
\\ 
\nonumber 
 & = & 
 E \left[ 
 \left( 
 \int_X h_1 ( \tau^{(k_{1,n})} ( x , \omega ) ) \omega ( dx ) 
 \right)^{l_1} 
 \cdots 
 \left( 
 \int_X h_m ( \tau^{(k_{m,n})} ( x , \omega ) ) \omega ( dx )  
 \right)^{l_m} 
 \right], 
\end{eqnarray} 
 converge to 
$$ 
 E \left[ 
 \left( 
 \int_X h_1 ( x ) \omega ( dx ) 
 \right)^{l_1} 
 \right] 
 \cdots 
 E \left[ 
 \left( 
 \int_X h_m ( x ) \omega ( dx )  
 \right)^{l_m} 
 \right] 
$$ 
 as $n$ goes to infinity. 
\\ 
 
 By Proposition~\ref{pr11} and the relation 
\begin{eqnarray} 
\nonumber 
\lefteqn{ 
\! \! \! \! \! \! \! \! \! \! \! \! \! \! \! \! \! \! \! \! \! \! \! 
\! \! \! \! \! \! \! \! \! \! \! \! \! \! \! \! \! \! \! 
 \varepsilon_{x_1,\ldots ,x_k}^+ 
 ( 
 u_1 (x_1 , \omega )  \cdots 
 u_k (x_k , \omega ) 
 ) 
 = 
 ( I + D_{x_1} ) \cdots ( I + D_{x_k} ) 
 ( 
 u_1 (x_1 , \omega )  \cdots 
 u_k (x_k , \omega ) 
 ) 
} 
\\ 
\label{djklddd} 
 & = & 
 \sum_{\Theta \subset \{1,\ldots ,k\} } 
 D_\Theta 
 ( 
 u_1 (x_1 , \omega )  \cdots 
 u_k (x_k , \omega ) 
 ), 
\end{eqnarray} 
 where 
 $D_{\Theta} = D_{x_1} \cdots D_{x_l}$ when 
 $\Theta = \{ x_1, \ldots , x_l\}$, 
 we can express the joint moment \eqref{djda} 
 as a finite sum of terms of the form 
\begin{equation} 
\label{fs} 
 E \left[ 
 \int_{X^k} 
 D_\Theta 
 \left( 
 \prod_{i_1\in Q_1} 
 h^{l^N_{1,i_1}}_{i_1} ( \tau^{(k_{i_1,n})} ( x_1 , \omega ) ) 
 \cdots 
 \prod_{i_k \in Q_k} 
 h^{l^N_{k,i_k}}_{i_k} ( \tau^{(k_{i_k,n})} ( x_k , \omega ) ) 
 \right) 
 \sigma (dx_1) 
 \cdots 
 \sigma (dx_k) 
 \right] 
, 
\end{equation} 
 where, with $N = l_1 + \cdots + l_m$,  
$l^N_{j,i} : = 
 | P^N_j \cap ( l_1+\cdots + l_{i-1},l_1+\cdots + l_i ]| 
$ 
 and 
\begin{equation} 
\label{qj} 
 Q_j = \{ i \in \{ 1, \ldots , m \} 
 \ : \ 
 l^N_{j,i} \geq 1 \}, 
 \qquad 
 j = 1, \ldots , k, 
\end{equation} 
 and $\Theta \subset \{x_1,\ldots ,x_k\}$. 
 Note that 
 when $\Theta = \{ x_1 , \ldots , x_k \}$, 
 Condition~\eqref{111.111.111} shows the vanishing of 
 \eqref{fs} due to the relation 
\begin{equation}
\label{eq:DDD0}
 D_{x_1} \cdots D_{x_k} \big( u_1( x_1 , \omega ) 
 \cdots u_k ( x_k , \omega )\big) =
 \sum_{\Theta_1\cup\cdots\cup \Theta_k = \{1, \ldots, k\}} D_{\Theta_1} 
 u_1( x_1 , \omega ) \cdots D_{\Theta_k} u_k ( x_k , \omega ),
\end{equation} 
 where the above sum includes all (possibly empty) sets 
 $\Theta_1, \ldots , \Theta_k$ whose union is 
 $\{x_1, \ldots, x_k\}$. 
 Hence in the sequel we can assume that 
 $\Theta = \{ x_1,\ldots , x_l \} \subset \{x_1,\ldots ,x_{k-1} \}$, 
\\ 
 
 The proof is split in four steps that 
 are based on the evaluation of \eqref{fs}. 
\\ 
 
\noindent 
{\em Step~1.} The term \eqref{fs} vanishes as $n$ tends to 
 infinity if $|Q_k| \geq 2$. 
\bigskip 
 
 If $Q_k$ contains at least two distinct indexes $a,b$ with 
 $1 \leq a < b \leq m$ and $l<k$, we have 
\begin{eqnarray} 
\nonumber 
\lefteqn{ 
\! \! \! \! \! \! \! \! \! \! \! \! \! \! \! \! 
 \int_X 
 \prod_{j \in Q_k} 
 h_j 
 ( \tau^{(k_{j,n})} ( x_k , {\omega} ) ) 
 \sigma (dx_k) 
 \leq 
 \int_X 
 h_a 
 ( \tau^{(k_{a,n})} ( x_k , {\omega} ) ) 
 h_b 
 ( \tau^{(k_{b,n})} ( x_k , {\omega} ) ) 
 \sigma (dx_k) 
} 
\\ 
\nonumber 
 & = & 
 \int_X 
 h_a 
 ( \tau^{(k_{a,n}-1)} ( \tau ( x_k , {\omega} ) , \tau_* {\omega} ) ) 
 h_b ( \tau^{(k_{b,n}-1)} ( \tau ( x_k , {\omega} ) , \tau_* {\omega} ) ) 
 \sigma (dx_k) 
\\ 
\nonumber 
 & = & 
 \int_X 
 h_a 
 ( \tau^{(k_{a,n}-1)} ( x_k , \tau_* {\omega} ) ) 
 h_b ( \tau^{(k_{b,n}-1)} ( x_k , \tau_* {\omega} ) ) 
 \sigma (dx_k) 
\\ 
\label{dkkdda} 
 & = & 
 \int_X 
 h_a 
 ( x_k ) 
 h_b ( \tau^{(k_{b,n}-k_{a,n})} ( x_k , \tau_*^{k_{a,n}} {\omega} ) ) 
 \sigma (dx_k) 
, 
\end{eqnarray} 
 $\omega \in \Omega$, 
 where used the invariance of $\sigma$ under 
 $\tau ( \cdot , \omega ) : X \longrightarrow X$. 
 By \eqref{dkkdda}, this shows that for all $p\geq 1$ we 
 have 
\begin{eqnarray} 
\nonumber 
\lefteqn{ 
 E\left[ 
 \left( 
 \int_X 
 \prod_{j \in Q_k} 
 h_j ( \tau^{(k_{j,n})} ( x_k , {\omega} ) ) 
 \sigma (dx_k) 
 \right)^p 
 \right] 
} 
\\ 
\nonumber 
 & \leq & 
 E\left[ 
 \left( 
 \int_X 
 h_a 
 ( x_k ) 
 h_b ( \tau^{(k_{b,n}-k_{a,n})} ( x_k , \tau_*^{k_{a,n}} {\omega} ) ) 
 \sigma (dx_k) 
 \right)^p 
 \right] 
\\ 
\nonumber 
 & = & 
 E\left[ 
 \left( 
 \int_X 
 h_a 
 ( x_k ) 
 h_b ( \tau^{(k_{b,n}-k_{a,n})} ( x_k , \omega ) ) 
 \sigma (dx_k) 
 \right)^p 
 \right] 
, 
\end{eqnarray} 
 while 
$$ 
 \int_X 
 h_a 
 ( x_k ) 
 h_b ( \tau^{(k_{b,n}-k_{a,n})} ( x_k , \omega ) ) 
 \sigma (dx_k) 
$$ 
 is a.s. bounded by 
 $\int_X 
 h_a 
 ( x_k ) 
 \sigma (dx_k)$ 
 and tends to zero in probability by \eqref{condfg} 
 as $n$ goes to infinity since $h_a$ has compact support and 
 $\lim_{n \to \infty} k_{b,n} - k_{a,n} = +\infty$. 
 Hence 
\begin{equation} 
\label{dkldd} 
 \lim_{n \to \infty} \int_X 
 \prod_{j \in Q_r} 
 h_j ( \tau^{(k_{j,n})} ( x_r , {\omega} ) ) 
 \sigma (dx_r) = 0 
\end{equation} 
 in $L^p (\Omega )$, $p\geq 1$. 
\\ 
 
 From \eqref{dkldd} and the fact 
 that $\Theta = \{ x_1,\ldots , x_l \} \subset \{x_1,\ldots ,x_{k-1} \}$ 
 it is apparent that 
 \eqref{fs} will tend to zero as $n$ tends to infinity, 
 however to conclude Step~1 we need to an integrability 
 argument. 
\\ 
 
 For this, using the relation 
$$ 
 D_\Theta 
 = 
 \sum_{\eta \subset \Theta } 
 (-1)^{|\eta |+l} 
 \epsilon^+_\eta, 
$$ 
 where $\epsilon^+_\eta$ is defined as in \eqref{djkld1}, 
 we rewrite \eqref{fs} as a linear combination of terms of the form 
$$ 
 E \left[ 
 \int_{X^k} 
 \varepsilon^+_\eta 
 \left( 
 \prod_{i_1\in Q_1} 
 h^{l^N_{1,i_1}}_{i_1} ( \tau^{(k_{i_1,n})} ( x_1 , \omega ) ) 
 \cdots 
 \prod_{i_k \in Q_k} 
 h^{l^N_{k,i_k}}_{i_k} ( \tau^{(k_{i_k,n})} ( x_k , \omega ) ) 
 \right) 
 \sigma (dx_1) 
 \cdots 
 \sigma (dx_k) 
 \right] 
, 
$$ 
 with $\eta = \{ x_1,\ldots , x_l \} \subset \{x_1,\ldots ,x_{k-1} \}$. 
 Applying the 
 first moment Mecke identity \eqref{mecke1} 
 to the variable $x_1$, we get 
\begin{eqnarray*} 
\lefteqn{ 
 E \left[ 
 \int_{X^k} 
 \varepsilon^+_\eta 
 \prod_{j=1}^k 
 \left( 
 \prod_{i_j\in Q_j} 
 h^{l^N_{j,i_j}}_{i_j} ( \tau^{(k_{i_j,n})} ( x_j , \omega ) ) 
 \right) 
 \sigma (dx_1) 
 \cdots 
 \sigma (dx_k) 
 \right] 
} 
\\ 
 & = & 
 E \left[ 
 \int_{X^{k-1}} \int_X  
 \varepsilon^+_{\eta \setminus \{x_1\} } 
 \left( 
 \prod_{j=1}^k 
 \left( 
 \prod_{i_j\in Q_j} 
 h^{l^N_{j,i_j}}_{i_j} ( \tau^{(k_{i_j,n})} ( x_j , \omega ) ) 
 \right) 
 \right) 
 \omega (dx_1) 
 \sigma (dx_2) 
 \cdots 
 \sigma (dx_k) 
 \right] 
\\ 
 & = & 
 E \left[ 
 \int_{X^{k-1}} 
 \varepsilon^+_{\eta \setminus \{ x_1 \} } 
 \left( 
 \int_X 
 \prod_{j=1}^k 
 \left( 
 \prod_{i_j\in Q_j} 
 h^{l^N_{j,i_j}}_{i_j} ( \tau^{(k_{i_j,n})} ( x_j , \omega ) ) 
 \right) 
 \omega (dx_1) 
 \right) 
 \sigma (dx_2) 
 \cdots 
 \sigma (dx_k) 
 \right] 
\\ 
 & & 
 \! \! \! \! \! \! \! \! \! \! \! \! \! \! \! \! \! \! \! \! \! 
 - 
 \sum_{r=2}^l 
 E \left[ 
 \int_{X^{k-1}} 
 \varepsilon^+_{\eta \setminus \{ x_1 \} } 
 \left( 
 \prod_{i_1\in Q_1} 
 h^{l^N_{1,i_1}}_{i_1} ( \tau^{(k_{i_1,n})} ( x_l , \omega ) ) 
 \prod_{j=2}^k 
 \left( 
 \prod_{i_j\in Q_j} 
 h^{l^N_{j,i_j}}_{i_j} ( \tau^{(k_{i_j,n})} ( x_j , \omega ) ) 
 \right) 
 \right) 
 \sigma (dx_2) 
 \cdots 
 \sigma (dx_k) 
 \right] 
, 
\end{eqnarray*} 
 where we used the relation 
$$ 
 \varepsilon^+_{x_2} \cdots \varepsilon^+_{x_l} 
 \int_X v(x_1 , \omega ) \omega (dx_1) 
 = 
 \int_X \varepsilon^+_{x_2} \cdots \varepsilon^+_{x_l} 
 v(x_1, \omega ) \omega (dx_1) 
 + 
 \sum_{r \in \eta \setminus \{ x_1 \} } 
 \varepsilon^+_{x_2} \cdots \varepsilon^+_{x_l} v(x_r, \omega ). 
$$ 
 After inductively exhausting all elements of $\eta$ 
 by repeating the above argument 
 we find that \eqref{fs} rewrites as a linear combination of 
 terms of the form 
\begin{equation} 
\label{dlddd} 
 E \left[ 
 \left( 
 \prod_{j=1}^{l'} 
 \int_X 
 \prod_{i_j\in R_j} 
 h^{l^N_{j,i_j}}_{i_j} ( \tau^{(k_{i_j,n})} ( x_j , \omega ) ) 
 \omega (dx_j ) 
 \right) 
 \left( 
 \prod_{j={l'}+1}^{k'} 
 \int_X 
 \prod_{i_j\in R_j} 
 h^{l^N_{j,i_j}}_{i_j} ( \tau^{(k_{i_j,n})} ( x_j , \omega ) ) 
 \sigma (dx_j ) 
 \right) 
 \right] 
, 
\end{equation} 
 $1 \leq {l'} < {k'}$, where 
 $\{R_1,\ldots ,R_{k'} \}$ is another family of subsets of 
 $\{ 1 , \ldots , m \}$ with $R_{k'} = Q_k$. 
\\ 
 
 Denoting by $K \subset X$ a compact set containing the supports 
 of $h_1,\ldots , h_m$, all $l$ terms in the 
 left product in \eqref{dlddd} are a.s. bounded by 
 the random variable  
\begin{eqnarray} 
\nonumber 
 \int_X {\bf 1}_K ( \tau^{(k_{i_j,n})} ( x_j , \omega ) ) 
 \omega ( dx_j ) 
 & = & 
 \int_X {\bf 1}_K ( 
 \tau^{(k_{i_j,n}-1)} ( \tau ( x_j , \omega ) , \tau_* \omega ) ) 
 \omega ( dx ) 
\\ 
\label{djld11} 
 & = & 
 \int_X {\bf 1}_K ( 
 \tau^{(k_{i_j,n}-1)} ( x_j , \tau_* \omega ) ) 
 \tau_* \omega ( dx ) 
, 
\end{eqnarray} 
 which has the same distribution as 
$\displaystyle 
 \int_X {\bf 1}_K ( 
 \tau^{(k_{i_j,n}-1)} ( x , \omega ) ) \omega ( dx ) 
$ 
 since $\tau_* : \Omega \longrightarrow \Omega$ leaves the 
 Poisson measures $\pi_\sigma$ invariant. 
 by Theorem~3.3 of \cite{prinv} or \cite{prigirid}. 
 By decreasing induction on 
 $k_{i_j,n},k_{i_j,n}-1,\ldots ,1$, this shows that \eqref{djld11} 
 has the Poisson distribution of 
$\displaystyle 
 \int_X {\bf 1}_K ( x ) \omega ( dx ) 
 = 
 \omega (K) 
$ 
 with parameter $\sigma (K) < \infty$, 
 in particular it has finite moments of all orders. 
\\ 
 
 On the other hand, 
 the terms of index $j = l+1,\ldots , {k'}-1$ in the right product 
 \eqref{dlddd} are uniformly bounded in $n$ by $\sigma (K)$ as 
 in 
 \eqref{dkkdda}, 
 and the last term of index 
 ${k'}$ converges to $0$ in $L^p(\Omega )$ for all $p \geq 1$ 
 by \eqref{dkldd} since $Q_{k}$ is not a singleton. 
 Hence by H\"older's inequality, \eqref{fs} tends to $0$ as $n$ 
 goes to infinity. 
\\ 
 
\noindent 
{\em Step~2.} 
 As a consequence of Step~1 we only 
 need to consider terms \eqref{fs} of the form 
\begin{eqnarray*} 
\lefteqn{ 
 E \left[ 
 \int_{X^k} 
 D_\Theta 
 \left( 
 h^{l^N_{1,i_1}}_{i_1} ( \tau^{(k_{i_1,n})} ( x_1 , \omega ) ) 
 \cdots 
 h^{l^N_{k,i_k}}_{i_k} ( \tau^{(k_{i_k,n})} ( x_k , \omega ) ) 
 \right) 
 \sigma (dx_1) 
 \cdots 
 \sigma (dx_k) 
 \right] 
} 
\\ 
 & = & 
 E \left[ 
 \int_{X^k} 
 D_\Theta 
 \left( 
 h^{l^N_{1,i_1}}_{i_1} ( \tau^{(k_{i_1,n})} ( x_1 , \omega ) ) 
 \cdots 
 h^{l^N_{k-1,i_{k-1}}}_{i_{k-1}} ( \tau^{(k_{i_{k-1},n})} ( x_{k-1} , \omega ) ) 
 \right) 
\right. 
\\ 
 & & 
 \ \ \ \ \ \ \ \ \ \ \ \ \ \ \ \ \ \ \ \ \ \ \ \ \ \  \ \ \ \ \ \ \ \ \ \ \ \ \  \ \ \ \ \ \ \ \ \ \ \ \ \ 
 \times h^{l^N_{k,i_k}}_{i_k} ( \tau^{(k_{i_k,n})} ( x_k , \omega ) ) 
 \sigma (dx_1) 
 \cdots 
 \sigma (dx_k) 
 \Biggr] 
\end{eqnarray*} 
 where 
 $Q_k = \{i_k\}$ is a singleton. 
 By invariance of $\tau ( \cdot , \tilde{\omega} ) : X \longrightarrow X$ 
 for any 
 $\tilde{\omega} \subset \omega \cup \{ x_1, \ldots , x_l \} 
 \in \Omega$, we have 
\begin{eqnarray} 
\nonumber 
\lefteqn{ 
 \! \! \! \! \! \! \! \! \! \! \! \! \! \! \! 
 \int_{X} 
 h^{l^N_{k,i_k}}_{i_k} ( \tau^{(k_{i_k,n})} ( x_k , \tilde{\omega} ) ) 
 \sigma ( dx_k ) 
 = 
 \int_{X} 
 h^{l^N_{k,i_k}}_{i_k} ( \tau^{(k_{i_k,n}-1)} ( \tau ( x_k , \tilde{\omega} ) \tau_* \tilde{\omega} ) )  
 \sigma ( dx_k ) 
} 
\\ 
\nonumber 
 & = & 
 \int_{X} 
 h^{l^N_{k,i_k}}_{i_k} ( \tau^{(k_{i_k,n}-1)} ( x_k , \tau_* \tilde{\omega} ) ) 
 \sigma ( dx_k ) 
 = 
 \int_{X} 
 h^{l^N_{k,i_k}}_{i_k} ( \tau ( x_k , \tau_*^{k_{i_k,n}-1} \tilde{\omega} ) ) 
 \sigma ( dx_k ) 
\\ 
\label{dhjklddw} 
 & = & \int_{X} h^{l^N_{k,i_k}}_{i_k} ( x_k ) \sigma ( dx_k ), 
\end{eqnarray} 
 where the step before last is reached by induction 
 on $1, \ldots , k_{i_k,n}-1$. 
 Since \eqref{dhjklddw} is deterministic, 
 the integral in 
 $\sigma (dx_k)$ can then be factored out of $D_\Theta$ in \eqref{fs} 
 and we can reconsider \eqref{fs} at the order $k-1$ instead of $k$. 
\\ 
 
\noindent 
{\em Step~3.} Decreasing induction on $k$. 
 
 After implementing Step~2, 
 from \eqref{eq:DDD0} and Condition~\eqref{111.111.111} 
 we can again assume that 
 $\Theta = \{ x_1,\ldots , x_l \} \subset \{x_1,\ldots ,x_{k-2} \}$, 
 and repeating Step~2 above by further decrementing $k$ 
 we find that \eqref{fs} vanishes as $n$ tends to infinity 
 unless $Q_j$ is a singleton for all $j=1,\ldots , k=m$ 
 and $\Theta$ is empty, in which case we have 
\begin{eqnarray*} 
\lefteqn{ 
 \! \! \! \! \! \! \! \! \! \! \! \! \! \! \! \! \! \! \! \! \! \! \! \! \! \! \! \! \! \! \! \! \! \! \! \! \! 
 \lim_{n\to \infty} 
 E \left[ 
 \int_{X^k} 
 \left( 
 h^{l^N_{1,i_1}}_{i_1} ( \tau^{(k_{i_1,n})} ( x_1 , \omega ) ) 
 \cdots 
 h^{l^N_{k,i_k}}_{i_k} ( \tau^{(k_{i_k,n})} ( x_k , \omega ) ) 
 \right) 
 \sigma (dx_1) 
 \cdots 
 \sigma (dx_k) 
 \right] 
} 
\\ 
\nonumber 
 & = & 
 \int_{X} 
 h_1^{l^N_{1,i_1}} (x ) 
 \sigma (dx) 
 \cdots 
 \int_X 
 h_{k_m}^{l^N_{m,i_m}} (x ) 
 \sigma (dx) 
. 
\\ 
\end{eqnarray*} 
\noindent 
{\em Step~4.} 
 To conclude, taking again $N = l_1 + \cdots + l_m$ we let 
$$ 
 U^i_j : = P^N_j \cap ( l_1+\cdots + l_{i-1},l_1+\cdots + l_i ], 
 \qquad 
 i = 1,\ldots , m, 
 \quad 
 j=1,\ldots ,k, 
$$ 
 and note that from \eqref{qj} and Step~3, 
 \eqref{fs} vanishes as $n$ tends to infinity, 
 unless $\Theta = \emptyset$ and the cardinal 
$$ 
 l^N_{j,i} = |U^i_j| = | P^N_j \cap ( l_1+\cdots + l_{i-1},l_1+\cdots + l_i ]| 
$$ 
 is either $0$ or $1$ 
 for all $j=1,\ldots ,k$ and $i\in \{ 1,\ldots , m\}$. 
\\ 
 
 Hence we only need to consider partitions of $\{1,\ldots , k\}$ 
 of the form 
$$ 
 \{ U^1_1, \ldots , U^1_{k_1},\ldots ,U^m_1,\ldots ,U^m_{k_m} \}
$$ 
 such that for all $i = 1,\ldots , m$, 
$$ 
 \{ U^i_1, \ldots , U^i_{k_i} \} 
$$ 
 is a partition of $\{l_{1}+\cdots + l_{i-1} +1, \ldots , l_1 + \cdots l_i \}$ 
 having (say) $k_i$ non empty sets, 
 $k_i \in \{1,\ldots,l_i \}$, 
 after a suitable re-indexing of the 
 lower index $j$ in $U^i_j$. 
\\ 
 
 Then by \eqref{djklddd} we have 
\begin{eqnarray*} 
\lefteqn{ 
 \lim_{n\to \infty} 
 E \left[ 
\left( 
 \int_X h_1 ( \tau^{(k_{1,n})} ( x , \omega ) ) \omega ( dx ) 
 \right)^{l_1} 
 \cdots 
 \left( 
 \int_X h_m ( \tau^{(k_{m,n})} ( x , \omega ) ) \omega ( dx )  
 \right)^{l_m} 
 \right] 
} 
\\ 
\nonumber 
 & = & 
 \lim_{n\to \infty} 
 \sum_{k=1}^N 
 \sum_{P^N_1,\ldots , P^N_k } 
 E \left[ 
 \int_{X^k} 
 \varepsilon_{x_1,\ldots ,x_k}^+ 
 \left( 
 \prod_{j=1}^k 
 \prod_{i=1}^m 
 h_i^{l^N_{i,j}} 
 ( \tau^{(k_{i,n})} ( x_j , \omega ) ) 
 \right) 
 \sigma (dx_1) \cdots  \sigma (dx_k) 
 \right] 
\\ 
\nonumber 
 & = & 
 \lim_{n\to \infty} 
 \sum_{k_1=1}^{l_1} 
 \cdots 
 \sum_{k_m=1}^{l_m} 
 \sum_{U^1_1 \cup \ldots \cup U^1_{k_1} = \{ 1 , \ldots , l_1 \}} 
 \cdots 
 \sum_{U^m_1 \cup \ldots \cup U^m_{k_m} = \{ l_1+\cdots +l_{m-1}+1 , \ldots , l_1+\cdots + l_m \}} 
\\ 
\nonumber 
& & 
 \! \! \! \! \! \! \! \! \! \! \! \! \! \! 
 E \left[ 
 \int_{X^{k_1+\cdots + k_m}} 
 \varepsilon^+_{x_1,\ldots ,x_{k_1+\cdots +k_m }} 
 \left( 
 \prod_{i=1}^m 
 \prod_{q_i=1}^{k_i} 
 h_i^{| U^i_{q_i} |} 
 ( \tau^{(k_{i,n})} ( x_{k_1+\cdots + k_{i-1} + q_i} , \omega ) ) 
 \right)  
 \sigma (dx_1) \cdots  \sigma (dx_{k_1+\cdots + k_m}) 
 \right] 
\\ 
\nonumber 
 & = & 
 \lim_{n\to \infty} 
 \sum_{k_1=1}^{l_1} 
 \cdots 
 \sum_{k_m=1}^{l_m} 
 \sum_{U^1_1 \cup \ldots \cup U^1_{k_1} = \{ 1 , \ldots , l_1 \}} 
 \cdots 
 \sum_{U^m_1 \cup \ldots \cup U^m_{k_m} = \{ l_1+\cdots +l_{m-1}+1 , \ldots , l_1+\cdots + l_m \}} 
 \sum_{\Theta \subset \{ 1,\ldots , k_1+\cdots +k_m \}} 
\\ 
\nonumber 
& & 
 E \left[ 
 \int_{X^{k_1+\cdots + k_m}} 
 D_\Theta 
 \left( 
 \prod_{i=1}^m 
 \prod_{q_i=1}^{k_i} 
 h_i^{| U^i_{q_i} |} 
 ( \tau^{(k_{i,n})} ( x_{k_1+\cdots + k_{i-1} + q_i} , \omega ) ) 
 \right)  
 \sigma (dx_1) \cdots  \sigma (dx_{k_1+\cdots + k_m}) 
 \right] 
\\ 
\nonumber 
 & = & 
 \sum_{k_1=1}^{l_1} 
 \cdots 
 \sum_{k_m=1}^{l_m} 
 \sum_{U^1_1 \cup \ldots \cup U^1_{k_1} = \{ 1 , \ldots , l_1 \}} 
 \cdots 
 \sum_{U^m_1 \cup \ldots \cup U^m_{k_m} = \{ 1 , \ldots , l_m \}} 
\\ 
 & & 
 \int_{X} 
 h_1^{|U^1_1|} (x ) 
 \sigma (dx)
 \cdots 
 \int_{X} 
 h_{k_1}^{|U^1_{k_1}|} (x ) 
 \sigma (dx)
 \cdots 
 \int_X 
 h_m^{|U^m_1|} (x ) 
 \sigma (dx) 
 \cdots 
 \int_X 
 h_{k_m}^{|U^m_{k_m}|} (x ) 
 \sigma (dx) 
\\ 
 & = & 
 E \left[ 
 \left( 
 \int_X h_1 ( x ) \omega ( dx ) 
 \right)^{l_1} 
 \right] 
 \cdots 
 E \left[ 
 \left( 
 \int_X h_m ( x ) \omega ( dx )  
 \right)^{l_m} 
 \right] 
, 
\end{eqnarray*} 
 showing 
 that $\tau_*$ is mixing of all 
 orders $n\geq 1$, by density in $L^2 ( \Omega , \pi_\sigma )$ 
 of the polynomials in $\int_X h(x) \omega (dx)$, 
 $h \in {\cal C}_c (X)$. 
\end{Proof} 
\section{Examples} 
\label{examples} 
 We consider a family of 
 examples satisfying the hypotheses 
 of Theorem~\ref{prnt}, based on transformations conditioned 
 by a random boundary. 
 We let $X = \real^d$ with norm $\Vert \cdot \Vert$ 
 and for all $\omega \in \Omega$ we denote by 
 $\omega_e \subset \omega$ denote the extremal vertices of 
 the convex hull of $\omega \cap B(0,1)$. 
 We also denote by ${\cal C} (\omega )$ the convex hull of $\omega$, 
 with interior $\dot{\cal C} (\omega )$. 
\\ 
 
 Consider a mapping $\widehat{\tau} : X \times \Omega \longrightarrow X$ 
 such that for all $\omega \in \Omega$, 
 $\widehat{\tau} ( \cdot , \omega ) : X \longrightarrow X$ leaves 
 $X \setminus \dot{{\cal C}}(\omega_e )$ invariant 
 (including the extremal vertices $\omega_e$ of ${\cal C} (\omega_e )$) 
 while 
 $\widehat{\tau} : \dot{{\cal C}}(\omega_e ) \times \Omega \longrightarrow \dot{{\cal C}}(\omega_e )$ 
 shifts the points inside $\dot{{\cal C}}(\omega_e )$ 
 depending on the data of $\omega_e$, i.e. we have 
\begin{equation} 
\label{adfdsfg} 
 \widehat{\tau} ( x , \omega ) = 
 \left\{ 
 \begin{array}{ll} 
 \widehat{\tau} ( x , \omega_e ), 
 & 
 x\in \dot{{\cal C}}(\omega_e ), 
\\ 
\\ 
 x , 
 & 
 x\in X\setminus \dot{{\cal C}}(\omega_e ). 
\end{array} 
\right. 
\end{equation} 
 As shown in Proposition~\ref{ffjkl} below, such a transformation 
 $\widehat{\tau}$ satisfies the vanishing condition \eqref{cyclic4} hence by 
 Theorem~3.3 of \cite{prinv} or \cite{prigirid} 
 the mapping 
 $\widehat{\tau}_* : \Omega \longrightarrow \Omega$ 
 leaves $\pi_\sigma$ invariant. 
 The next figure shows an example of behaviour 
 such a transformation, with a 
 finite set of points for simplicity of illustration. 
\\ 
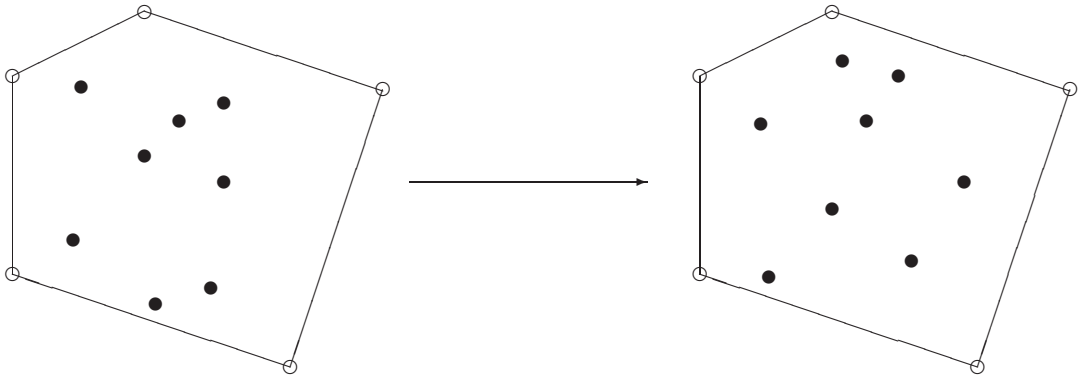
\vspace{0.4cm} 
\begin{center} 
\begin{picture}(230,130)(-60,-90)
\linethickness{0.2pt}
\put(-140,40){\circle{5}} 
\put(-140,-35){\circle{5}} 
\put(-140,-35){\line(0,1){75}} 
\put(-140,40){\line(2,1){50}} 
\put(-90,64.5){\circle{5}} 
\put(-90,64.5){\line(3,-1){90}} 
\put(0,35){\circle{5}} 
\put(0,35){\line(-1,-3){35}} 
\put(-35,-70){\circle{5}} 
\put(-35,-70){\line(-3,1){105}} 
\put(-60,0){\circle*{5}} 
\put(-65,-40){\circle*{5}} 
\put(-60,30){\circle*{5}} 
\put(-90,10){\circle*{5}} 
\put(-117,-22){\circle*{5}} 
\put(-77,23){\circle*{5}} 
\put(-86,-46){\circle*{5}} 
\put(-114,36){\circle*{5}} 

\put(10,0){\vector(3,0){90}} 

\put(120,40){\circle{5}} 
\put(120,-35){\circle{5}} 
\put(120,-35){\line(0,1){75}} 
\put(120,40){\line(2,1){50}} 
\put(170,64.5){\circle{5}} 
\put(170,64.5){\line(3,-1){90}} 
\put(260,35){\circle{5}} 
\put(260,35){\line(-1,-3){35}} 
\put(225,-70){\circle{5}} 
\put(225,-70){\line(-3,1){105}} 
\put(220,0){\circle*{5}} 
\put(195,40){\circle*{5}} 
\put(200,-30){\circle*{5}} 
\put(170,-10){\circle*{5}} 
\put(143,22){\circle*{5}} 
\put(183,23){\circle*{5}} 
\put(174,46){\circle*{5}} 
\put(146,-36){\circle*{5}} 
\end{picture}
\end{center} 
\vspace{-0.5cm} 
 
 Using the mapping $\widehat{\tau} : X \times \Omega \longrightarrow X$, 
 we will build examples of interacting transformations 
 $\tau : X \times \Omega \longrightarrow X$ that satisfy 
 Conditions~\eqref{111.111.111} and \eqref{condfg}. 
\subsubsection*{Vanishing condition \eqref{111.111.111}} 
\begin{prop} 
\label{ffjkl} 
 Let $\widehat{\tau} : X \times \Omega \longrightarrow X$ 
 satisfy \eqref{adfdsfg} 
 and let 
 $f :X\longrightarrow X$ 
 be a bijective deterministic mapping 
 that preserves set convexity. 
 Then the transformation 
\begin{eqnarray} 
\nonumber 
 \tau & : X \times \Omega & \longrightarrow \ \Omega  
\\ 
\label{flllff} 
 & (\omega , x ) & \longmapsto \ 
 \tau ( x , \omega ) 
 : = f ( \widehat{\tau} ( x , \omega ) ) 
\end{eqnarray} 
 satisfies the vanishing condition \eqref{111.111.111}. 
\end{prop} 
\begin{Proof} 
 In order to check that \eqref{111.111.111} 
 holds for all $m\geq 1$, we note that 
 by induction on $k\geq 1$ we have 
\begin{equation} 
\label{jklf} 
 \tau^{(k)} ( x , \omega ) = 
 \tau^{(k)} ( x , \omega_e ), 
 \qquad 
 x \in X, 
\end{equation} 
 i.e. $\tau^{(k)} ( x , \omega )$ depends only on $x$ and on the points in 
 $\omega_e$. 
 Indeed, Relation~\eqref{jklf} is satisfied for $k=1$ by \eqref{adfdsfg} 
 and we have 
$$ 
 \tau^{(k+1)} ( x , \omega ) 
 = 
 \tau^{(k)} ( \tau ( x , \omega ) , \tau_* \omega ) 
 = 
 \tau^{(k)} ( \tau ( x , \omega_e ) , ( \tau_* \omega )_e ), 
$$ 
 while the positions of the points in $( \tau_* \omega )_e$ 
 themselves depend only on $\omega_e$ through the function $f$, 
 showing that $\tau^{(k+1)} ( x , \omega )$ depends only on 
 $\omega_e$ and $x$. 
\\ 
 
 On the other hand we can also show by induction that 
\begin{equation} 
\label{abi} 
 \tau^{(k)} ( x , \omega ) = f^k ( x ), \qquad 
 x \in X\setminus {\cal C} ( \omega_e ), 
\end{equation} 
 Indeed this condition is satisfied for $k=1$ by \eqref{adfdsfg} 
 and \eqref{flllff}. 
 Now since $f :X\longrightarrow X$ preserves 
 set convexity we have 
$$ 
 {\cal C} ( ( \tau_* \omega )_e ) 
 = 
 {\cal C} ( f ( \omega_e ) ) 
 \subset 
 f ( {\cal C} ( \omega_e ) ) 
, 
$$ 
 because $f ( {\cal C} ( \omega_e ) )$ 
 is convex and contains $f ( \omega_e )$, 
 hence since $f$ is bijective we get 
$$ 
 \tau ( x , \omega ) \in {\cal C} ( ( \tau_* \omega )_e ) 
 \Longrightarrow 
 {\tau} ( x , \omega ) \in f ( {\cal C} ( \omega_e ) ) 
 \Longrightarrow 
 \widehat{\tau} ( x , \omega ) \in {\cal C} ( \omega_e ) 
 \Longrightarrow 
 x \in {\cal C} ( \omega_e ) 
, 
$$ 
 i.e. 
\begin{equation} 
\label{abc} 
 x \in X \setminus {\cal C} ( \omega_e ) \Longrightarrow 
 \tau ( x , \omega ) = f(x) \in X\setminus {\cal C} ( ( \tau_* \omega )_e ), 
 \qquad 
 x \in X. 
\end{equation} 
 Therefore, assuming that \eqref{abi} holds at the rank $n\geq 1$, 
 for every $x \in X \setminus {\cal C} ( \omega_e )$ we get, 
 by \eqref{abc}, 
$$ 
 \tau^{(k+1)} ( x , \omega ) 
 = 
 \tau^{(k)} ( \tau ( x , \omega ) , \tau_* \omega ) 
 = 
 f^k ( \tau ( x , \omega ) ) 
 = 
 f^{k+1} ( x ), 
$$ 
 which is \eqref{abi} at the rank $k+1$. 
 In the remainder of this proof we will conclude 
 from \eqref{jklf} and \eqref{abi} as in Proposition~3.3 
 of \cite{bretonprivaultfact} 
 and \cite{prinv} 
 that the vanishing Condition~\eqref{111.111.111} is satisfied, 
 i.e. we show that 
\begin{equation} 
\label{djdkd1.0} 
 D_{\Theta_1} 
 \tau^{(k_1)} ( x_1 , \omega ) 
 \cdots 
 D_{\Theta_m} 
 \tau^{(k_m)} ( x_m , \omega ) 
 = 
 0, 
\end{equation} 
 for every family $\{ \Theta_1, \ldots , \Theta_m \}$ 
 of (non empty) subsets such that 
 $\Theta_1 \cup \cdots \cup \Theta_m = \{x_1,\ldots ,x_m \}$, 
 $x_1,\ldots ,x_m \in X$ and all 
 $\pi_\sigma ( d \omega )$-a.s., 
 $k_1,\ldots ,k_m \geq 1$,  $m\geq 1$. 
\\ 
 
 Note that whenever 
 $x_i$ lies inside of ${\cal C}(\omega)={\cal C}(\omega_e)$ 
 then by \eqref{jklf} we have 
\begin{eqnarray*}
D_{x_i} \tau^{(k)} ( x_j , \omega)
&=&\tau^{(k)} (x_j,\omega\cup\{x_i\})-\tau^{(k)} (x_j,\omega)
=\tau^{(k)} (x_j,(\omega\cup\{x_i\})_e)-\tau^{(k)} (x_j,\omega_e)\\
&=&\tau^{(k)} (x_j,\omega_e)-\tau^{(k)} (x_j,\omega_e)= 0 
\end{eqnarray*}
 for all $i,j=1,\ldots ,m$ and $k\geq 1$, 
 hence 
$ 
D_\eta \tau^{(k)} ( x_j , \omega)
 = 
0
$ 
 provided $ \{x_i \} \subset \eta \subset \{ x_1,\ldots , x_m \}$. 
\\ 
 
 Consequently it suffices to consider the case 
 where ${\cal C}(\omega \cup \{ x_1,\ldots ,x_m \})$ 
 has (at least) one extremal point denoted $x_e$ 
 within $\{ x_1,\ldots ,x_m \}$. 
\\ 
 
 Now, for all $\eta \subset \{ x_1,\ldots , x_m\}$ we have 
$$ 
 \tau^{(k)} ( x_e , \omega \cup \eta ) = 
 \tau^{(k)} ( x_e , \omega ) = f^k (x_e) 
$$ 
 by 
 \eqref{abi}, hence 
$$ 
 D_\Theta \tau^{(k)} ( x_e , \omega ) = 0, 
$$ 
 for all $\Theta \subset \{ x_1,\ldots , x_m\}$, 
 due to the relation 
\begin{eqnarray*} 
 D_\Theta 
 \tau^{(k)} ( x_e , \omega ) 
 & = & 
 \sum_{\eta \subset \Theta } 
 (-1)^{ | \Theta | + 1 - | \eta |} 
 \tau^{(k)} ( x_e , \omega \cup \eta ) 
\\ 
 & = & 
 f^k ( x_e ) 
 \sum_{\eta \subset \Theta } 
 (-1)^{ | \Theta | + 1 - | \eta |} 
\\ 
 & = & 
 f^k ( x_e ) 
 ( 1 - 1 )^{|\Theta |+1} 
\\ 
 & = & 0, 
\end{eqnarray*} 
 where the summation above holds over all (possibly empty) 
 subset $\eta$ of $\Theta$. 
 As a consequence, a factor in \eqref{djdkd1.0} 
 has to vanish. 
\end{Proof} 
\subsubsection*{Zero-type condition \eqref{condfg}} 
 In order for the zero-type condition \eqref{condfg}  
 to hold it suffices that 
$$ 
 \lim_{n\to \infty} 
 \Vert \tau^{(n)} ( x , \omega ) \Vert 
 = \infty, \qquad 
 \omega \in \Omega, \quad x\in \real^d. 
$$  
 For this 
 we can assume for example that 
 $\tau : X \times \Omega \longrightarrow X$ 
 satisfies a random dilation property 
\begin{equation} 
\label{gjlgjl} 
 \Vert 
 \tau ( x , \omega ) 
 \Vert \geq C (\omega ) \Vert x \Vert_s, 
 \qquad 
 \omega \in \Omega, \quad x \in \real^d, 
\end{equation} 
 for a random variable 
$$ 
 C : \Omega \longrightarrow (1,\infty )
. 
$$ 
 In this case, for any $g , h\in {\cal C}_c (X)$ with support 
 in $B(0,r)$ for some $r>0$, we have 
$$ 
 \lim_{n\to \infty} 
 \langle g , h \circ \tau^{(n)} \rangle_{L^2_\sigma (X)} = 0 
, 
 \qquad 
 \omega \in \Omega, 
$$ 
 because the support of $x \longmapsto h ( \tau^{(n)} ( x , \omega ) )$ 
 is in $B(0, r C^{-n} (\omega ) )$ by construction, for all 
 $\omega \in \Omega$. 
\\ 
 
 Condition~\eqref{gjlgjl} holds in particular when 
 $f:\real^d \longrightarrow \real^d$ in \eqref{flllff} 
 satisfies the dilation property 
$$ 
 \Vert f(x) \Vert \geq r \Vert x \Vert, 
 \qquad 
 x \in \real^d, 
$$ 
 for some $r > 1$, and 
 $\widehat{\tau} : X \times \Omega \longrightarrow X$ satisfies 
$$ 
 \Vert \widehat{\tau} ( x , \omega ) \Vert \geq c ( \omega ) 
 \Vert x \Vert, 
 \qquad 
 \omega \in \Omega,
 \quad 
 x \in \real^d, 
$$ 
 for some $r > 1$ and $c : \Omega \longrightarrow (0,1]$ such that 
 $\inf_{\omega \in \Omega} c ( \omega ) > 1/r$. 
\\ 
 
 For example in case $f(x) = r U x$, $x\in \real^d$, 
 where $r>1$ and $U:\real^d \longrightarrow \real^d$ is a linear isometry 
 of $\real^d$, 
 the intensity measure $\sigma ( dx) : = \Vert x \Vert^{-d} dx$ 
 is invariant by $f : \real^d \longrightarrow \real^d$, and 
 if $c(\omega )=1$, 
 the measure-preserving mapping 
 $\widehat{\tau} ( \cdot , \omega_e ) : X \longrightarrow X$ can be built 
 from any isometric transformations of $\dot{{\cal C}}( \omega_e )$. 
 This includes for example any random rotation 
 within a (random) disk contained in $\dot{{\cal C}}( \omega_e )$. 
 
\footnotesize 

\def\cprime{$'$} \def\polhk#1{\setbox0=\hbox{#1}{\ooalign{\hidewidth
  \lower1.5ex\hbox{`}\hidewidth\crcr\unhbox0}}}
  \def\polhk#1{\setbox0=\hbox{#1}{\ooalign{\hidewidth
  \lower1.5ex\hbox{`}\hidewidth\crcr\unhbox0}}} \def\cprime{$'$}

\end{document}